\documentclass[a4paper, 12pt]{amsart}

\usepackage{amssymb}
\usepackage{amsthm}
\usepackage{amsmath}
\usepackage[mathscr]{eucal}
\usepackage{comment}

\AtBeginDocument{%
   \def\MR#1{}
}

\theoremstyle{plain}
\newtheorem{thm}{Theorem}[section]		
\newtheorem{prop}[thm]{Proposition}
\newtheorem{cor}[thm]{Corollary}
\newtheorem{lem}[thm]{Lemma}

\theoremstyle{definition}
\newtheorem{df}{Definition}[section]

\theoremstyle{remark}
\newtheorem{rmk}{Remark}[section]
\newtheorem*{ac}{Acknowledgements}

\newcommand{\zz}{\mathbb{Z}}
\newcommand{\qq}{\mathbb{Q}}
\newcommand{\rr}{\mathbb{R}}

\DeclareMathOperator{\card}{Card}

\DeclareMathOperator{\met}{Met}

\DeclareMathOperator{\ult}{UMet}

\DeclareMathOperator{\metdis}{\mathcal{D}}

\DeclareMathOperator{\umetdis}{\mathcal{UD}}

\newcommand{\myens}{\mathbb{E}}

\newcommand{\yosub}{\subset}

\newcommand{\yogset}{\mathrm{Neb}}

\makeatletter
\@addtoreset{equation}{section}

\makeatother

\newcommand{\yodc}[1]{\mathrm{DC}(#1)}
\newcommand{\yogl}[1]{\mathrm{GL}(#1)}

\DeclareMathOperator{\yodiam}{diam}

\newcommand{\yoceil}[1]{\lceil #1 \rceil}

\newcommand{\yozset}{\mathcal{Z}}

\newcommand{\yoclassf}{\mathscr{F}}

\newcommand{\yoclasst}{\mathscr{T}}

\newcommand{\yoclassc}{\mathscr{C}}

\newcommand{\yoclassq}{\mathscr{Q}}

\newcommand{\yopcantor}{\Lambda}

\newcommand{\younivset}[1]{\mathrm{T}(#1)}
\newcommand{\younivsetc}[1]{\mathrm{CT}(#1)}

\newcommand{\yocalu}{\mathscr{N}}

\newcommand{\myou}{\mathbb{O}}

\begin{document}

\title[On comeager sets]
{
On comeager sets of metrics whose ranges are disconnected
}

\author[Yoshito Ishiki]
{Yoshito Ishiki}
\address[Yoshito Ishiki]
{\endgraf
Photonics Control Technology Team
\endgraf
RIKEN Center for Advanced Photonics
\endgraf
2-1 Hirasawa, Wako, Saitama 351-0198, Japan}
\email{yoshito.ishiki@riken.jp}

\date{\today}

\maketitle

\begin{abstract}
For a metrizable space $X$, 
we denote by 
$\mathrm{Met}(X)$
the space of all metric that generate the same topology of 
$X$. 
The space $\mathrm{Met}(X)$ is 
equipped with the supremum distance. 
In this paper, 
for every strongly zero-dimensional 
metrizable space $X$, 
we prove that 
the set of all metrics whose ranges are  
 closed totally disconnected subsets of the line 
is a dense  $G_{\delta}$ 
subspace 
in 
$\mathrm{Met}(X)$. 
As its application, 
we show that some 
sets of universal metrics are meager in 
spaces of metrics. 
\end{abstract}

\section{Introduction}

A metric $d$ on $X$ is said to be 
\emph{ultrametric} if it satisfies 
$d(x, y)\le d(x, z)\lor d(z, y)$ for all $x, y, z\in X$, 
where $\lor$ is the maximum operator on $\rr$. 
A topological space is said to be 
\emph{metrizable} (resp.~\emph{ultrametrizable})
 if there exists a metric (resp.~ultrametric)
that generates the same topology of the space. 
Let $X$ be a metrizable space. 
Let $S$ be a subset of $[0, \infty)$ with 
$0\in S$. 
We denote by $\met(X; S)$
(resp.~$\ult(X; S)$)
the set of all metrics (resp.~ultrametrics) that 
generate the same topology of $X$ taking 
values in $S$. 
We often write 
$\met(X)=\met(X; [0, \infty))$. 
We define a map 
$\metdis_{X}\colon \met(X)^{2}\to [0, \infty]$ by 
$\metdis_{X}(d, e)=\sup|d(x, y)-e(x, y)|$.
Then $\metdis_{X}$  is a metric on $\met(X)$ taking values in $[0, \infty]$. 
As in the case of ordinary 
metric spaces, using open 
balls, 
we can introduce the topology on
$(\met(X), \metdis_{X})$. 
We can also define an ultrametric 
$\umetdis_{X}^{S}$ on $\ult(X; S)$. 
We omit its definition since we do not use it in this paper.

A topological space 
$X$ is said to be \emph{strongly $0$-dimensional}
if for every pair $A, B$ of disjoint closed subsets of $X$, 
there exists a clopen set $V$ such that $A\yosub V$
and $V\cap B=\emptyset$. 
Such a space is  sometimes said to be \emph{ultranormal}.
Note that 
a topological space $X$ is 
ultrametrizable if and only if 
it is metrizable and strongly $0$-dimensional
(see \cite[Theorem II]{MR80905}). 

In \cite{Ishiki2020int}, \cite{ishiki2021dense}
and \cite{Ishiki2021ultra}, 
the author investigated 
geometric properties in the space 
$(\met(X), \metdis_{X})$ of metrics
and the space
 $(\ult(X; S), \umetdis_{X}^{S})$ of 
 ultrametrics. 

There are some researches concerning the range of a metric and 
the topology of an underlying set.  
Dovgoshey--Shcherbak
\cite{MR4335845}
proved 
that 
an ultrametrizable space $X$ is 
separable if and only if 
for every $d\in \ult(X; \rr_{\ge 0})$, 
the set $\{\, d(x, y)\mid x, y\in X\, \}$
is countable. 

Broughan \cite{MR314012}
proved that the following 
statements 
are equivalent to each other:
\begin{enumerate}
\item The space $X$ is ultrametrizable;
\item 
$\met(X; H)\neq \emptyset$, 
where  $H=\{0\}\cup \{\, 1/n\mid n\in \zz_{\ge 1}\, \}$; 
\item 
There exists a decreasing sequence 
$\{a_{i}\}_{i\in \zz_{\ge 1}}$ converging to $0$, 
such that 
$\met(X; B)\neq \emptyset$, 
where $B=\{0\}\cup \{\, a_{i}\mid i\in \zz_{\ge 1}\, \}$. 
\end{enumerate}

A metric $d$
on $X$ is said to be \emph{gap-like} if 
for every $p\in X$, 
the set $\{\, d(p, x)\,\mid x\in X \}$ is 
not dense in any neighborhood of $0$ in 
$[0, \infty)$. 
Broughan 
\cite[Theorem 7]{MR651427} proved that 
the Euclidean metric on $\rr\setminus \qq$
is a uniform limit of
gap-like  metrics on $\rr\setminus \qq$. 
We improve Broughan's result on  approximation of a metric by gap-like metrics. 

A topological space is said to be 
\emph{totally disconnected} if 
every  its connected component is a singleton. 
We denote by 
$\yozset$
the set of all closed totally disconnected
subsets of $[0, \infty)$ containing $0$. 
We define $\yodc{X}=\bigcup_{S\in \yozset}\met(X; S)$. 
Namely, $\yodc{X}$ is the set of all 
metrics whose ranges are closed and totally  disconnected. 
The following is our first result:
\begin{thm}\label{thm:main}
Let $X$ be a strongly 
$0$-dimensional 
metrizable space. 
Then the set $\yodc{X}$ is 
dense $G_{\delta}$ in the space 
$(\met(X), \metdis_{X})$. 
\end{thm}

A subset of topological space is 
said to be 
\emph{nowhere dense} if 
its closure has no interior points. 
A subset of a topological space is 
\emph{meager} if it is a countable  union of 
 nowhere dense subsets of this space. 
A subset of a topological space is 
said to be \emph{comeager} if its 
complement is meager. 
Note that if a subset of a topological space $X$
 contains a dense $G_{\delta}$ set, 
 then it is comeager. If $X$ is a Baire space, 
 then the converse is true. 

Let $X$ be a metrizable space. 
We denote by 
$\yogl{X}$ the set of all 
gap-like metrics in $\met(X)$. 
As a consequence of Theorem \ref{thm:main}, 
we have:
\begin{thm}\label{thm:main2}
Let $X$ be a strongly 
$0$-dimensional 
metrizable space. 
Then the set $\yogl{X}$ is 
comeager in the space 
$(\met(X), \metdis_{X})$. 
\end{thm}
Theorem 
 \ref{thm:main}
and \ref{thm:main2} will be proven in 
Section \ref{sec:proof}. 
Our  main results
(Theorems \ref{thm:main}
and \ref{thm:main2})
can be considered as  generalizations of 
Broughan's theorem 
\cite[Theorem 7]{MR651427}. 
\begin{rmk}
The author does not know whether 
the set $\yodc{X}$ has 
the anti-transmissible property defined in \cite{Ishiki2020int} or not. 
\end{rmk}

In Section \ref{sec:applications}, 
as applications of Theorem 
 \ref{thm:main}, 
we  show 
that some
sets of  universal metrics are meager in 
space of metrics. 

\section{Proofs of Theorems}\label{sec:proof}

For a metric space $(X, d)$ and a subset $A$, 
we denote by $\yodiam_{d}(A)$ the 
diameter of $A$ with respect to $d$. 
We begin with an amalgamation of metrics. 

\begin{prop}\label{prop:amalfi}
Let $I$ be a set. 
Let $(X, d)$ be a metric space. 
Let $\{B_i\}_{i\in I}$ be a covering  of $X$  
consisting of  mutually disjoint  clopen  subsets 
of $X$. 
Let $P=\{p_i\}_{i\in I}$ be points with $p_i\in B_i$. 
Let $\{e_i\}_{i\in I}$ be 
a set of metrics such that 
$e_i\in \met(B_i)$. 
Let $h$ be a metric on $P$ 
generating  the discrete topology on $P$.  
We define a function $D:X^2\to [0, \infty)$ by 
\[
D(x, y)=
\begin{cases}
e_i(x, y) & \text{ if $x, y\in B_i$;}\\
e_i(x, p_i)+h(p_i, p_j)+e_j(p_j, y)  & \text{if $x\in B_i$ and $y\in B_j$.}
\end{cases}
\]
Then 
$D\in \met(X)$ 
and 
$D|_{B_i^2}=e_i$ 
for all $i\in I$. 
Moreover, 
if for every $i\in I$ 
we have 
$\yodiam_{d}(B_i)\le \epsilon$
 and 
 $\yodiam_{e_i}(B_i)\le \epsilon$, 
 then 
$\metdis_{X}(D, d)\le 4\epsilon+\metdis_{P}(d|_{P^2}, h)$. 
\end{prop}
\begin{proof}
The statement  that $D$ is a metric  can be proven in a similar way to 
\cite[Proposition 3.1]{ishiki2021dense}. 
We now prove that $D$ generates the same topology 
of $X$.  
Take
$x\in X$. Then there  exists a unique element  $i\in I$ with 
 $x\in B_{i}$. 
 Since $h$ generates the  discrete topology on $P$, 
there exists $\delta_{i}\in (0, \infty)$ such that 
$\delta_{i}\le h(p_{i}, p_{j})$ for all $j\in I$ with $i \neq j$. 

Take $r\in (0, \infty)$ with $r<\delta_{i}$.
By the definition of $D$, 
 the open ball centered at $x$ with radius $r$
with respect to $D$ is contained in $B_{i}$. 
Since $X$ is homeomorphic to the direct sum of 
$\{B_{i}\}_{i\in I}$, we conclude that $D\in \met(X)$. 

To prove the latter part, 
we take $x, y\in X$. 
If $x, y\in B_{i}$ for some $i\in I$, then, by $\yodiam_{e_{i}}(B_{i})\le \epsilon$
and $\yodiam_{d}(B_{i})\le \epsilon$, 
we have 
\[
|D(x, y)-d(x, y)|=|e_{i}(x, y)-d(x, y)|\le 2\epsilon<
4\epsilon+\metdis_{P}(d|_{P^2}, h).
\] 
If $x\in B_{i}$ and $y\in B_{j}$ for some distinct $i, j\in I$, then 
we have 
\begin{align*}
|D(x, y)-d(x, y)|
&\le e(x, p_{i})+e(y, p_{j})+|h(p_{i}, p_{j})-d(x, y)|\\
&\le 
2\epsilon +|h(p_{i}, p_{j})-d(x, y)|. 
\end{align*}
We also have 
\begin{align*}
&|h(p_{i}, p_{j})-d(x, y)| \le
|d(x, y)-d(p_{i}, p_{j})|+|h(p_{i}, p_{j})-d(p_{i}, p_{j})|\\
&\le d(x, p_i)+d(x, p_{j})+\metdis_{P}(d|_{P^2}, h)
\le 2\epsilon +\metdis_{P}(d|_{P^2}, h). 
\end{align*}
This implies that $|D(x, y)-d(x, y)|\le 4\epsilon+\metdis_{P}(d|_{P^2}, h)$. 
Therefore we conclude that $\metdis_{X}(D, d)\le 4\epsilon+\metdis_{P}(d|_{P^2}, h)$.
\end{proof}

The following is deduced from 
\cite[Proposition 1.2 and Corollary 1.4]{ellis1970extending}. 
\begin{lem}\label{lem:ultraparacompact}
Let $X$ be a strongly $0$-dimensional 
metrizable space. 
Then every open cover $\{U_{i}\}_{i\in I}$ of $X$ has a
refinement $\{E_{j}\}_{j\in J}$ such that 
each $E_{j}$ is open in $X$, and 
$\bigcup_{j\in J}E_{j}=X$,  and 
$E_{j}\cap E_{j^{\prime}}=\emptyset$ if 
$j\neq j^{\prime}$. 
In this case, each $E_{j}$ is clopen. 
\end{lem}

For each $x\in \rr$, 
we denote by $\yoceil{x}$
the minimum integer of all $k\in \zz$
with $x\le k$. 
The following two lemmas are related to 
\emph{metric-preserving functions}
(see \cite{dobovs1998metric}). 
\begin{lem}\label{lem:subadditive}
For all $x, y\in \rr$, 
we have 
$\yoceil{x+y}\le \yoceil{x}+\yoceil{y}$. 
\end{lem}
\begin{proof}
By $x\le \yoceil{x}$ and 
$y\le \yoceil{y}$, 
we have $x+y\le \yoceil{x}+\yoceil{y}$. 
By $\yoceil{x}+\yoceil{y}\in \zz$ and 
by the definition of $\yoceil{x+y}$, 
we conclude that 
$\yoceil{x+y}\le \yoceil{x}+\yoceil{y}$. 
\end{proof}

The following is 
a reformulation of 
\cite[Theorem 4.1]{MR3027678}
or
\cite[Theorem 1 and Proposition 1]{dobovs1998metric}. 
\begin{thm}\label{thm:subadditive}
Let $f\colon [0, \infty)\to [0, \infty)$ be an increasing function such that 
$f^{-1}(0)=0$. 
Then the following statements are equivalent to each other:
\begin{enumerate}
\item For every  set $X$ and for every metric $d$ on $X$, 
the function   $f\circ d\colon X\times X\to [0, \infty)$ is 
a metric on $X$.
\item The function $f$ is subadditive, i.e., 
the inequality  $f(x+y)\le f(x)+f(y)$ holds 
for all $x, y\in [0, \infty)$. 
\end{enumerate}
\end{thm}

\begin{rmk}
\cite[Theorem 4.1]{MR3027678} is 
proven using  only  the property that $\rr$ is a linearly ordered 
Abelian group. 
Thus, 
Theorem \ref{thm:subadditive} is still true 
for metrics taking values in linearly ordered Abelian groups. 
For more discussion of  generalized  metrics, 
we refer the readers to \cite{ishiki2022sim}. 
\end{rmk}

For $\eta\in (0, \infty)$, 
we define 
$\eta\cdot \zz=\{\, \eta\cdot n\mid n\in \zz\, \}$ and 
$\eta\cdot \zz_{\ge 0}=(\eta\cdot \zz)\cap \rr_{\ge 0}$. 
\begin{prop}\label{prop:discretemetric}
Let $X$ be a discrete topological space, 
$\eta\in (0, \infty)$ 
and $d\in \met(X)$.
 Then, there exists 
a metric $e\in \met(X)$ 
such that 
 $\metdis_{X}(d, e)\le \eta$
 and 
 $e(x, y)\in \eta\cdot \zz_{\ge 0}$ 
 and $\eta\le e(x, y)$
 for all distinct 
 $x, y\in X$. 
\end{prop}
\begin{proof}
We
define $e\colon X\times X\to \rr$ by 
$e(x, y)=\eta\cdot 
\yoceil{\eta^{-1}\cdot d(x, y)}$ if 
$x\neq y$; otherwise $e(x, x)=0$. 
According to Theorem \ref{thm:subadditive}, 
the map $e$ is a metric on $X$. 
Since $d$ generates the  discrete topology on $X$, 
by the definition of $\yoceil{x}$, 
we conclude that $\eta\le e(x, y)$
 for all distinct $x, y\in X$. 
 Thus,  $e\in \met(X)$ is a  required metric. 
\end{proof}

Let $A$ and $B$ be subsets of $\rr$. 
We define the sum of $A$ and $B$ by 
$A+B=\{\, a+b\mid a\in A, b\in B\, \}$. 
\begin{df}\label{df:eeeset}

Let 
$\eta\in (0, \infty)$ and 
$u\in (0, 1)$. 
Let us define 
the sets 
$\myou(\eta, u)$
and $\myens(\eta, u)$
as 
\[
\myou(\eta, u)=\{0\}\cup \{\eta\cdot u^{n}\mid n\in \zz_{\ge 0}\}
\]
and
\begin{align*}
\myens(\eta, u)=
\eta\cdot \zz_{\ge 0}+
\myou(\eta, u)+\myou(\eta, u). 
\end{align*}
Note that the set $\myens(\eta, u)$ is 
the closure of 
the set 
$ \{\, 
\eta\cdot (l+u^{n}+u^{m})\mid l, n, m\in \zz_{\ge 0}\, \}$.  
\end{df}

\begin{lem}\label{lem:enszzz}
Let 
$\eta\in \rr_{\ge 0}$ and 
$u\in (0, 1)$. Then 
we have $\myens(\eta, u)\in \yozset$. 
\end{lem}
\begin{proof}
Since
$\myens(\eta, u)$ is countable, 
it is totally disconnected. 
Since $\myens(\eta, u)$ is 
the sum of the
 compact set $\myou(\eta, u)+\myou(\eta, u)$ and
the   closed set $\eta\cdot \zz_{\ge 0}$, 
it is closed (see 
\cite[Corollary 1, $\S$ 4.1, p.251]{MR1726779} or 
\cite[Theorem 1.4.30]{MR2433295}). 
Thus $\myens(\eta, u)\in \yozset$. 
\end{proof}

\begin{thm}\label{thm:dense}
Let $X$ be a 
strongly $0$-dimensional metrizable space. 
Let $d\in \met(X)$. 
Let $\epsilon\in (0, \infty)$. 
Then there exist
$r\in (0, 1)$ and 
 a metric $D$ on $X$
such that: 
\begin{enumerate}
\item\label{item:inens}
$D\in \met(X; \myens(\epsilon/5, r))$; 
\item\label{item:epsilon}
$\metdis_{X}(d, D)\le \epsilon$. 
\end{enumerate}
In particular, 
the set $\yodc{X}$ is dense in 
$\met(X)$. 
\end{thm}
\begin{proof}
We put $\eta=\epsilon/5$ and $r=\min\{1/2, \epsilon/10\}$. 
Applying Lemma \ref{lem:ultraparacompact}
to a covering of $X$ consisting of open balls 
with radius $r$, 
we can find 
a clopen covering $\{E_{i}\}_{i\in I}$ of $X$. 
For each $i\in I$, 
we
 take a point $p_{i}\in E_{i}$. 
  In this case, 
 we have $\yodiam_{d}(E_{i})\le 2r\le \eta$. 
 Put $Y=\{p_{i}\}_{i\in I}$. 
By Proposition \ref{prop:discretemetric},  
there exists
 $h\in \met(Y; \eta\cdot \zz_{\ge 0})$
 such that $h$ generates the discrete topology on $Y$ and $\metdis_{Y}(d|_{Y^{2}}, h)\le \eta$. 
 Put $R=\{0\}\cup \{\, \eta\cdot r^{n}\mid n\in \zz_{\ge 0}\, \} (=\myou(\eta, r))$. 
According to \cite[Proposition 2.14]{Ishiki2021ultra}, 
we can take a metric $e_{i}\in \ult(E_{i}; R)$. 
Then we have  $\yodiam_{e_{i}}(E_{i})\le \eta$. 
We define 
\[
D(x, y)=
\begin{cases}
e_i(x, y) & \text{if $x, y\in E_{i}$;}\\
e_i(x, p_i)+h(p_i, p_j)+e_j(p_j, y)  & \text{if $x\in E_{i}$ and $y\in E_{j}$.}
\end{cases}
\]
Applying  Proposition \ref{prop:amalfi} to 
$\{p_{i}\}_{i\in I}$, 
$\{E_{i}\}_{i\in I}$,
$\{e_{i}\}_{i\in I}$, 
 $h$,  and $\eta$, 
we have  $D \in \met(X)$ and 
we obtain 
$\metdis_{X}(D, d)\le 4\eta+\metdis_{Y}(d|_{Y^{2}}, h)\le 5\eta=\epsilon$. 
By the definition of $D$, we have $D(x, y)\in \myens(\eta, r)$ for all $x, y\in X$. 
This complete the proof of Theorem \ref{thm:dense}. 
\end{proof}

\begin{df}\label{df:nebula}
Let $q\in \zz_{\ge 0}$. 
We say that 
a 
closed subset $A$ of $[0, \infty)$
is 
a
\emph{$q$-nebula}
if
there exists a 
family
$\{I_{i}\}_{i=0}^{k}$
 of intervals 
satisfying the 
 the following conditions:
\begin{enumerate}
\item 
we have
$0\in I_{0}$ and 
$A=\bigcup_{i=1}^{k}I_{i}$;
\item 
each $I_{i}$  is 
a closed interval in $[0, \infty)$;
\item 
for all $i\in \{0, \dots, k-1\}$, 
we have 
$\yodiam(I_{i})< 2^{-q}$, where ``$\yodiam$'' stands for the 
diameter with respect to the Euclidean metric; 
\item 
the set $I_{k}$ is unbounded and 
$I_{k}\yosub (q, \infty)$; 
\item
if $i\neq j$, 
then we have $I_{i}\cap I_{j}=\emptyset$.
\end{enumerate}
\end{df}

\begin{lem}\label{lem:intersectionnebula}
Let $\{\, A_{q}\mid q\in \zz_{\ge 0}\, \}$
be a family of subsets of $[0, \infty)$ such that 
$A_{q}$ is  a $q$-nebula for every $q\in \zz_{\ge 0}$. 
 Then the set 
 $\bigcap_{q\in \zz_{\ge 0}}A_{q}$ 
 belongs to  $\yozset$. 
\end{lem}
\begin{proof}
Put $S=\bigcap_{q\in \zz_{\ge 0}}A_{q}$. Then $S$ is closed. 
Since  every  $q\in \zz_{\ge 0}$ satisfies $0\in A_{q}$, 
we have $0\in S$. 
For the sake of contradiction, we
assume that 
 $S$ has a 
connected component containing at least two points. 
Then there is $a, b\in S$
such that $[a, b]\yosub S$. 
Take $q\in \zz_{\ge 0}$
so that $b<q$ and $2^{-q}\le |b-a|$. 
Then $[a, b]\yosub A_{q}$. 
This is a contradiction to the definition of 
the $q$-nebula. 
Thus, 
$S$ is totally disconnected. 
Hence 
$S\in \yozset$. 
\end{proof}

\begin{lem}\label{lem:intervalzero}
Let $S$ be a subset of $\rr$. 
Let $a, b\in \rr$. 
If $S$ is totally disconnected  and 
 $a<b$, 
then 
$[a, b]\setminus S\neq \emptyset$. 
\end{lem}
\begin{proof}
If $[a, b]\yosub S$, then 
$S$ has a connected component containing at least two points.  
\end{proof}

\begin{lem}\label{lem:extqnebula}
Let $S\in \yozset$. 
Then 
for every $q\in \zz_{\ge 0}$, 
there exists a $q$-nebula $A$ such that $S\yosub A$ and  
each compact connected component of $A$ intersects $S$. 
\end{lem}
\begin{proof}
Let $M=(q+1)2^{q+1}$. 
For each  $m\in \{0, \dots, M\}$, 
we put $C_{m}=m\cdot 2^{-(q+1)}$. 
Put $\eta=2^{-(q+3)}$
Put $A_{m}=C_{m}-\eta$ and 
$B_{i}=C_{m}+\eta$. 
Put $t_{0}=0$. 
According to 
Lemma \ref{lem:intervalzero}, 
for each  $m\in \{1, \dots, M\}$, 
we can 
take $t_{m}\in [A_{m}, B_{m}]\setminus S$. 
Then, by the definitions  of $A_{m}$ $B_{m}$, and $\eta$, 
we obtain  $|t_m-t_{m+1}|<2^{-q}$ for all $m$
and $q<t_{M}$. 
Take $k\in \zz_{\ge 0}$
and a map $\phi\colon \{0, \dots, k\}\to \{0, \dots, M\}$
such that  $\phi(0)=0$ and  $\phi(k)=M$ and $[t_{\phi(i)}, t_{\phi(i)+1}]\cap S\neq \emptyset$
for all $i\in \{0, \dots, k-1\}$
and 
$S\yosub [t_{\phi(k)}, \infty)\cup  \bigcup_{i=0}^{k-1}[t_{\phi(i)}, t_{\phi(i)+1}]$. 
Since $S\cap [0, t_{\phi(k)}]$ is compact, 
for each $i\in \{0, \dots, k-1\}$, 
we can define $a_{i}$ and $b_{i}$ by 
$a_{i}=\min S\cap [t_{\phi(i)}, t_{\phi(i)+1}]$
and $b_{i}=\max S\cap [t_{\phi(i)}, t_{\phi(i)+1}]$. 
We define $a_{k}=\min S\cap [t_{\phi(k)}, \infty)$ if that set is non-empty; otherwise $a_{k}=t_{\phi(k)}$. 
Then
we have 
$[a_{0}, b_{0}]\yosub [0, t_{1})$ and 
 $[a_{i}, b_{i}]\yosub (t_{\phi(i)}, t_{\phi(i)+1})$ for all $i\in \{1, \dots, k-1\}$, 
and $[a_{k}, \infty)\yosub (q, \infty)$. 
For each $i\in\{0, \dots, k-1\}$, 
we put $I_{i}=[a_{i}, b_{i}]$
and $I_{k}=[a_{k}, \infty)$. 
Put $A=\bigcup_{i=0}^{k}I_{i}$. 
Then $A$ is a $q$-nebula with 
$S\yosub A$. 
For each $i\in \{0, \dots, k-1\}$, we have  $I_{i}\cap S\neq \emptyset$. 
Thus the set $A$
is a desired one. 
\end{proof}

\begin{prop}\label{prop:nebula}
Let $S\in \yozset$. Then there exists a family 
$\{A_{q}\}_{q\in \zz_{\ge 0}}$ such that 
each $A_{q}$ is a $q$-nebula 
and $S=\bigcap_{q\in \zz_{\ge 0}}A_{q}$. 
\end{prop}
\begin{proof}
According to Lemma \ref{lem:extqnebula}, 
for each $q\in \zz_{\ge 0}$, 
we can take a $q$-nebula $A_{q}$ such that $S\yosub A_{q}$
and each compact 
connected component intersects $S$.
Then we obtain 
$S\yosub \bigcap_{q\in \zz_{\ge 0}}A_{q}$.  
We shall show the converse inclusion. 
Let 
$x$ belong to $\bigcap_{q\in \zz_{\ge 0}}A_{q}$
and let 
$N\in \zz_{\ge 0}$ satisfy the inequality $N>x$. 
Then,  for every 
$i>N$, 
the point 
$x$ 
belongs to 
a compact connected component of 
$A_{i}$, 
say 
$C_{i}$. 
Note that 
$\yodiam C_{i}<2^{-i}$. 
By the assumption  on compact connected components of
$\{A_{q}\}_{q\in \zz_{\ge 0}}$, 
we have 
$S\cap C_{i}\neq \emptyset$ for all $i>N$. 
We can take 
$y_{i}\in S\cap C_{i}$. 
By 
$\yodiam C_{i}<2^{-i}$, we have 
$\lim_{i\to \infty}y_{i}=x$. 
Since $S$ is closed, we obtain 
$x\in S$. 
Therefore we conclude that 
$S=\bigcap_{q\in \zz_{\ge 0}}A_{q}$. 
\end{proof}

We omit the proof of the following lemma
since it is elementary. 
\begin{lem}\label{lem:elementary}
Let $a, b \in \rr$ and $\epsilon\in (0, \infty)$. 
If $x\in [a, b]$ and $|x-y|\le \epsilon$, 
then 
$y\in [a-\epsilon, b+\epsilon]$. 
\end{lem}
\begin{df}
Let $X$ be a metrizable space. 
We denote by 
$\yogset_{q}(X)$ the 
set of 
all metrics $d\in \met(X)$
such that there exists
a $q$-nebula   $A$
with 
$d\in \met(X; A)$. 
\end{df}
\begin{prop}\label{prop:open}
Let $X$ be a 
metrizable space. 
Then 
for every $q\in \zz_{\ge 0}$, 
the set 
 $\yogset_{q}(X)$ is open in $\met(X)$. 
\end{prop}
\begin{proof}
Take $d\in \yogset_{q}(X)$. 
Then there exists a $q$-nebula $A$ with 
$d\in \met(X; A)$. 
Let $\{I_{i}\}_{i=0}^{k}$ be a 
family of closed interval 
appearing in  the definition of nebulae such that 
$A=\bigcup_{i=0}^{k}I_{i}$. 
For $i \in \{1, \dots, k-1\}$, 
we put $I_{i}=[a_{i}, b_{i}]$, 
 and 
 $I_{0}=[0, b_{0}]$ and
 $I_{k}=[a_{k}, \infty)$.
We can take $c\in (0, \infty)$ such that 
 if $i\neq j$, then for all $x\in I_{i}$ and $y\in I_{j}$, 
 we have $c<d(x, y)$. 
 We take a sufficient small $\epsilon\in (0, \infty)$
 so that $|b_{i}-a_{i}|+2\epsilon <2^{-q}$
 for all $i\in \{0, \dots, k-1\}$, 
 and $\epsilon<a_{k}-q$ and 
 $\epsilon <c/4$. 
We define  $J_{0}=[0, b_{0}+\epsilon]$. 
For each $i\in \{1, \dots, k-1\}$, 
we define $J_{i}=[a_{i}-\epsilon, b_{i}+\epsilon]$. 
We define $J_{k}=[a_{k}-\epsilon, \infty)$. 
Then the set $B=\bigcup_{i=0}^{k}J_{i}$ is 
a $q$-nebula. 
Take $e\in \met(X)$ with  $\metdis_{X}(d, e)< \epsilon$. 
Then, 
by the definition of $\epsilon$ and 
by Lemma \ref{lem:elementary},
we have  $e\in \met(X; B)$, and hence 
 $e\in \yogset_{q}(X)$. 
This means that $\yogset_{q}(X)$ is open in $\met(X)$. 
\end{proof}

\begin{proof}[Proof of Theorem \ref{thm:main}]
Let $X$ be a strongly $0$-dimensional metrizable space. 
By Lemma \ref{lem:enszzz} and  Theorem \ref{thm:dense}, 
the set $\yodc{X}$ is dense in $\met(X)$. 
Put $L=\bigcap_{q\in \zz_{\ge 0}}\yogset_{q}(X)$. 
We shall prove that 
$\yodc{X}=L$. 
Take $d\in \yodc{X}$. 
Then there exists $S\in \yozset$ with $d\in \met(X; S)$. 
By Proposition \ref{prop:nebula}, there exists a sequence 
$\{A_{q}\}_{q\in \zz_{\ge 0}}$ such that each $A_{q}$ is a $q$-nebula 
and $\bigcap_{q\in \zz_{\ge 0}}A_{q}=S$. 
For each $q\in \zz_{\ge 0}$, we have $d\in \yogset_{q}(X)$, 
we obtain 
$d\in L$. 
Thus 
$\yodc{X}\yosub L$. 
To prove the converse inclusion, we take 
$d\in L$. 
Then there exists a sequence $\{A_{q}\}$ such that each $A_{q}$ is 
a $q$-nebula and $d\in \met(X; A_{q})$. 
Put $T=\bigcap_{q\in \zz_{\ge 0}}A_{q}$. 
By Lemma \ref{lem:intersectionnebula}, 
the set $T\in \yozset$. 
We also have $d\in \met(X; T)$, 
and hence $d\in \yodc{X}$. 
Therefore we conclude that 
$\yodc{X}=L$. 
Since Proposition \ref{prop:open} states that 
each $\yogset_{q}(X)$ is open in $\met(X)$, 
the set $\yodc{X}$ is $G_{\delta}$ in $\met(X)$. 
This finishes the proof of Theorem \ref{thm:main}. 
\end{proof}

\begin{proof}[Proof of Theorem \ref{thm:main2}]
Since each $S\in \yozset$ is not dense in any neighborhood of $0$, 
we have 
$\yodc{X}\yosub \yogl{X}$. 
This finishes the proof of Theorem \ref{thm:main2}. 
\end{proof}

In \cite[Lemma 5.1]{Ishiki2020int}, 
it is proven that 
for 
a second countable locally compact metrizable space $X$, the space $(\met(X), \metdis_{X})$ is a Baire space. 
Using this fact,  we obtain the following corollary. 
\begin{cor}
Let $X$ be a strongly $0$-dimensional 
$\sigma$-compact locally compact metrizable space. 
Then the sets $\yodc{X}$ and $\yogl{X}$
have 
the second category in $(\met(X), \metdis_{X})$ in the sense of Baire. 
\end{cor}


\section{Applications}\label{sec:applications}
In this section, 
as applications of 
Theorem \ref{thm:main}, 
we will show that 
some sets of  universal metrics are meager. 
\subsection{Constructions of universal metrics}
Let $\mathscr{M}$ be a class of metric spaces. 
A metric space $(X, d)$ is 
\emph{universal for 
$\mathscr{M}$},  
or 
\emph{$\mathscr{M}$-universal} if 
for every $(Y, e)\in \mathscr{M}$, 
there exists an isometric embedding 
$I\colon (Y, e)\to (X, d)$. 
In this case we will say that the metric 
$d\colon X\times X\to [0, \infty)$ also is 
$\mathscr{M}$-universal. 

A metric space is call a \emph{finite metric space} if 
its cardinality is finite. 
Let $\yoclassf$ be the class of all finite metric spaces.

We shall construct some universal metrics. 
The following theorem is 
proven 
 by 
Holsztynski \cite{Hol} (see also \cite{Ov}). 
\begin{thm}\label{thm:Hol}
There exists an $\yoclassf$-universal 
metric in $\met(\rr)$. 
\end{thm}
In this paper, we prove a generalization of 
Theorem \ref{thm:Hol}. 

Denote by $\Gamma$ be the Cantor set. 
We say that a topological space is 
a \emph{Cantor space} if it is homeomorphic to 
$\Gamma$. 
We say that a topological space is 
a 
\emph{punctured Cantor space}
 if it is homeomorphic to 
$\Gamma\setminus \{p\}$ for some  $p\in \Gamma$. 
Note that  a topological space is a punctured Cantor space if and only if  it is 
homeomorphic to the countable disjoint  union of 
Cantor spaces.

We  prove that a punctured Cantor space has 
an $\yoclassf$-universal metric. 

Let $r\in (0, \infty)$. 
A subset $S$ of a metric space $(X, d)$ is 
said to be 
\emph{$r$-separated} if 
$r\le d(x, y)$ for all 
distinct 
$x, y\in S$. 

The following 
lemmas 
(Lemmas \ref{lem:univmetn}
and \ref{lem:finebd})
 are proven in 
\cite{Ov}. 
For the sake of self-containedness, 
we provide proofs. 
\begin{lem}\label{lem:univmetn}
Let $r\in (0, \infty)$. 
Let $(X, d)$ and $(Y, e)$ be metric spaces. 
Let $f\colon X\to Y$ be a surjective continuous map.
We define a function 
$\rho\colon X\times X\to [0, \infty)$ by 
\[
\rho(x, y)=\min\{d(x, y), r\}\lor e(f(x), f(y)), 
\] 
where the symbol $\lor$ stands for 
the maximal operator on $\rr$. 
Then $\rho\in \met(X)$ and the metric space $(X, \rho)$ is universal for all $r$-separated finite subspace of $Y$. 
\end{lem}
\begin{proof}
Since $\min\{d(x, y), r\}\in \met(X)$ and $f$ is continuous, 
we have $\rho\in \met(X)$. 
We next prove the universality of $\rho$. 
Take an arbitrary $r$-separated subset $E=\{p_{i}\mid i\in \{0, \dots, k\}\}$ of $Y$. 
For each $i$, we take $t_{i}\in X$ with $f(t_{i})=p_{i}$. 
We define $I\colon (E, e|_{E^{2}})\to (X, \rho)$ by 
$I(p_{i})=t_{i}$. 
Since $E$ is $r$-separated, we have 
$r\le e(f(t_{i}), f(t_{j}))$ for distinct $i, j$. 
Then we have $\rho(I(p_{i}), I(p_{j}))=\rho(t_{i}, t_{j})=e(f(t_{i}), f(t_{j}))=e(p_{i}, p_{j})$ for all $i, j$. Thus $I$ is an isometric embedding. 
Therefore we conclude that $\rho$ is universal for 
all
 $r$-separated finite subspaces of $Y$. 
\end{proof}

Let $n\in \zz_{\ge 1}$. 
We denote by $\yoclassc_{n}$ the 
class of all
$(n^{-1})$-separated 
 finite metric spaces 
 $(X, d)$
satisfying that 
  $\yodiam_{d}(X)\le n$, 
and 
$\card(X)\le n$, where 
``$\card$'' stands for 
the cardinality.

\begin{lem}\label{lem:finebd}
Let $n\in \zz_{\ge 1}$. 
Then the space $[0, n]^{n}$
equipped with the $\ell^{\infty}$-Euclidean metric 
 is $\yoclassc_{n}$-universal. 
\end{lem}
\begin{proof}
Let $(X, d)$ be a metric space in 
$\yoclassc_{n}$. 
We represent $X=\{\, p_{i}\mid i\in \{1, \dots, k\}\, \}$, 
where $k\le n$. 
For each $i\in \{1, \dots, k\}$, we define 
$F_{i}\colon \{1, \dots, n\}\to [0, n]$ by 
$F_{i}(m)=d(p_{i}, p_{m})$ if $i\le k$; 
otherwise, $F_{i}=0$. 
We define  $\phi\colon X\to [0, n]^{n}$ 
by $\phi(p_{i})=(F_{1}(i), \dots, F_{n}(i))$. 
As is the case of the
Fr\'{e}chet  embedding, 
using 
the triangle inequality,  
we conclude that $\phi$ is an isometric embedding. 
\end{proof}

According to Lemmas 
\ref{lem:univmetn} and 
\ref{lem:finebd}, 
we obtain: 
\begin{cor}\label{cor:efuniv}
Let $X$ be a metrizable space. 
Let $n\in \zz_{\ge 1}$. 
If there exists a continuous surjective map 
$f\colon X\to [0, n]^{n}$, 
then there exists a $\yoclassc_{n}$-universal 
metric $e\in \met(X)$. 
\end{cor}

\begin{prop}\label{prop:funivpc}
Let $\yopcantor$ be a punctured 
Cantor space. 
There exists an $\yoclassf$-universal 
metric in 
$\met(\yopcantor)$. 
\end{prop}
\begin{proof}
Take a sequence $\{K_{i}\}_{i\in \zz_{\ge 0}}$ 
of subsets of $\yopcantor$
satisfying  that 
 $\yopcantor=\coprod_{i\in \zz_{\ge 0}}K_{i}$ and 
 each $K_{i}$ is a Cantor space. 
Since every compact metrizable space is a continuous image of a Cantor space (see, for example,  
\cite[Theorem 30.7]{W1970}), we can take a 
continuous surjective map 
$f_{i}\colon K_{i}\to [0, n]^{n}$. 
Then, by Corollary \ref{cor:efuniv}, we can take 
a $\yoclassc_{n}$-universal 
metric  $e_{i}\in \met(K_{i})$. 
For each $i\in \zz_{\ge 0}$, 
we take $p_{i}\in K_{i}$. 
We
define a metric $h$ on $\{p_{i}\}_{i\in \zz_{\ge 0}}$
by $h(p_{i}, p_{j})=1$ if $i\neq j$. 
Applying  Proposition \ref{prop:amalfi} to 
$\{p_{i}\}_{i\in \zz_{\ge 0}}$, 
$\{K_{i}\}_{i\in \zz_{\ge 0}}$, 
$\{e_{i}\}_{i\in \zz_{\ge 0}}$, and $h$,
we can take 
$D\in \met(\yopcantor)$ such that 
$D|_{K_{i}^{2}}=e_{i}$. 
Since every metric space  in $\yoclassf$ 
belongs to $\yoclassc_{n}$ for some $n\in \zz_{\ge 1}$, 
we conclude that  
the metric $D$ is $\yoclassf$-universal. 
\end{proof}

The next is Hausdorff's metric extension theorem \cite{Ha1930}
(see also \cite{MR321026}). 
\begin{thm}\label{thm:Hausdorff}
For a metrizable space $X$, 
and for a closed subset $A$ of $X$, 
and for every $d\in \met(A)$, 
there exists $D\in \met(X)$
such that $D|_{A^{2}}=d$. 
\end{thm}

Using  Proposition \ref{prop:funivpc} and 
Theorem \ref{thm:Hausdorff}
we obtain the following  generalization of Theorem \ref{thm:Hol}. 
\begin{cor}\label{cor:funiv}
If a metrizable  space $X$ contains a 
punctured Cantor space as a closed subset, 
then there exists an $\yoclassf$-universal 
 metric in $\met(X)$. 
\end{cor}

\begin{rmk}
If a metrizable topological space $X$ is
  Polish and there is an unbounded metric
  $d\in \met(X)$ such that the complement of every bounded subspace of $(X, d)$ has uncountable cardinality, then $X$  contains a 
punctured Cantor space as closed subset. 
This follows from the fact that every uncountable Polish space
contains a Cantor space 
(see \cite[Corollary 6.5]{MR1321597}). 
\end{rmk}

Let $S$ be a subset of $[0, \infty)$ with 
$0\in S$. 
We denote by $\yoclasst(S)$ the class of all 
two-point metric
 spaces
  whose metrics take values in 
$S$.  

Some examples and properties of 
$\yoclasst(S)$-universal
metric  spaces can be 
found in \cite{bilet2015minimal} and \cite{holurl}. 

\begin{prop}\label{prop:univtexist}
Let $S$ be a countable subset of 
$[0, \infty)$ with $0\in S$. 
Let $X$ be a countable discrete space. 
Then there exists a 
$\yoclasst(S)$-universal 
metric in  $\met(X)$. 
\end{prop}
\begin{proof}
Put $S=\{0\}\cup \{s_{i}\}_{i\in \zz_{\ge 0}}$. 
Take subsets $A, B$ of $X$ such that $A\cap B=\emptyset$ and 
$X=A\sqcup B$ and $A$ and $B$ are countable. 
Put $A=\{a_{i}\}_{i\in \zz_{\ge 0}}$ and 
$B=\{b_{i}\}_{i\in \zz_{\ge 0}}$. 
For each $i\in \zz_{\ge 0}$, 
we define a metric $e_{i}$ on $\{a_{i}, b_{i}\}$ by 
$e_{i}(a_{i}, b_{i})=s_{i}$. 
We define a metric $h$ on $A$ such that 
$h(x, y)=1$ if $x\neq y$. 
Applying Proposition
\ref{prop:amalfi}  to 
 $\{a_{i}\}_{i\in \zz_{\ge 0}}$,
$\{\{a_{i}, b_{i}\}\}_{i\in \zz_{\ge 0}}$, 
$\{e_{i}\}_{i\in \zz_{\ge 0}}$, 
and $h$, 
we obtain a metric $D\in \met(X)$ such that 
$D(a_{i}, b_{i})=e_{i}(a_{i}, b_{i})=s_{i}$. 
Then 
$D$ is a desired one. 
\end{proof}

Proposition \ref{prop:univtexist} and Theorem 
\ref{thm:Hausdorff} implies:
\begin{cor}\label{cor:univexists}
Let $S$ be a countable subset of 
$[0, \infty)$ with  $0\in S$. 
If a topological space $X$ contains 
a countable discrete space
 as a closed subset, 
 then there exists 
a $\yoclasst(S)$-universal 
metric in $\met(X)$. 
\end{cor}

\subsection{The meagerness of sets of  universal metrics}

Let $X$ be a metrizable space. 
Let $\younivset{X; S}$ be the set of all 
$\yoclasst(S)$-universal metrics in $\met(X)$. 
Let 
$\younivsetc{X; S}$ be the 
 closure of $\younivset{X; S}$ in 
$\met(X)$. 

\begin{lem}\label{lem:univdenserange}
Let $S$ be a dense subset of $[0, \infty)$
with  $0\in S$. 
Let $X$ be a 
metrizable space. 
Then 
for every  $d\in \younivsetc{X; S}$, 
the set $\{\, d(x, y)\mid x, y\in X\, \}$ is 
dense in $[0, \infty)$. 
\end{lem}
\begin{proof}
By the definition of $\younivsetc{X; S}$, 
for $d\in \younivsetc{X; S}$, 
we can take 
 a sequence
 $\{e_{n}\}_{n\in \zz_{\ge 0}}$ in 
 $\younivset{X; S}$
satisfying 
 that $\metdis_{X}(d, e_{n})\le 2^{-n}$ for all $n\in \zz_{\ge 0}$. 
Take arbitrary 
 $q\in [0, \infty)$ and $\epsilon \in (0, \infty)$. 
Take $r\in S$ with 
$|q-r|\le \epsilon/2$ and
take  $k\in \zz_{\ge 0}$
with $2^{-k}\le \epsilon/2$. 
Since $e_{k}$ is $\yoclasst(S)$-universal, 
we can
take $x, y\in X$ such that 
$e_{k}(x, y)=r$. 
Then we have 
$|d(x, y)-r|\le 2^{-k}$, 
and hence 
$|d(x, y)-q|\le \epsilon$. 
Thus the set 
$\{\, d(x, y)\mid x, y\in X\, \}$ is dense in 
$[0, \infty)$. 
\end{proof}

The following is an application of 
Theorem \ref{thm:main}. 
\begin{thm}\label{thm:allication1}
Let $S$ be a dense subset of $[0, \infty)$
with $0\in S$. 
Let $X$ be a strongly $0$-dimensional 
metrizable space. 
Then the set $\younivsetc{X; S}$ is meager in $\met(X)$. 
\end{thm}
\begin{proof}
Since every $S\in \yozset$ is not dense in $[0, \infty)$, 
Lemma \ref{lem:univdenserange} implies that 
$\younivsetc{X; S}\cap \yodc{X}=\emptyset$. 
Since $\yodc{X}$ is comeager 
(see Theorem \ref{thm:main}), 
we conclude that $\younivsetc{X; S}$ is meager. 
\end{proof}

Let $\mathscr{M}$ be a class of finite metric spaces. 
We say that a metric space $(X, d)$ is 
 \emph{$\mathscr{M}$-injective}
 if 
for every  metric space $(A, m)\in \mathscr{M}$
and for every $B\subset A$, 
every isometric embedding 
$\phi\colon (B, m|_{B^{2}})\to (X, d)$ can be 
extended to an isometric embedding 
$\Phi \colon (A, m)\to (X, d)$. 

Let $R$ be a subset of $[0, \infty)$. 
Let $\yocalu(R)$ be the class of 
all finite ultrametric spaces whose metrics take values   in $R$. 
Let $\yoclassq$ be the class of all finite metric spaces 
whose metrics take values in $\qq$. 

For each subset  $R$ of $[0, \infty)$ with $0\in R$, 
there exists 
a complete  $\yocalu(R)$-injective  ultrametric space. 
If $R$ is countable, then it is unique up to isometry,  
and it is called the  the 
\emph{$R$-Urysohn universal ultrametric space} (see \cite{MR2754373} and \cite{MR4282005}). 
A  countable 
$\yoclassq$-injective metric space 
 uniquely exists up to isometry, 
and it is called the  the 
\emph{rational Urysohn universal metric space}
(see, for example, \cite{MR2278687}). 
 \begin{rmk}
 The completion of the rational Urysohn universal metric spaces is a complete separable $\yoclassf$-injective metric space, and it  is called the 
 \emph{Urysohn universal metric space}. 
 For more discussions on this space, we refer the readers to,  for example,  \cite{MR2278687}, 
 \cite{MR2435145} and 
 \cite{MR2277969}. 
 \end{rmk}

In this paper, 
we use only the fact
 that 
 $\yocalu(R)$-injective ultrametric spaces
and 
the rational Urysohn universal space are
$\yoclasst(R)$-universal 
and $\yoclasst(\qq_{\ge 0})$-universal, respectively. 
By Theorem \ref{thm:allication1}, 
we obtain the following two corollaries. 
\begin{cor}\label{cor:univunivuniv}
Let $X$ be a strongly $0$-dimensional 
metrizable space. 
Let $R$ be a dense subset of 
$[0, \infty)$ with $0\in R$. 
Then the following subsets of $\met(X)$ are
meager in $\met(X)$:
\begin{enumerate}
\item 
The set of all $\yoclassf$-universal metrics in 
$\met(X)$. 
\item 
The set of all metrics $d\in \met(X)$ such that 
$(X, d)$ is an $\yocalu(R)$-injective 
ultrametric 
space. 
\end{enumerate}
\end{cor}
\begin{rmk}
It can happen that 
the two sets appearing in 
Corollary \ref{cor:univunivuniv} are
 empty. 
\end{rmk}

Note that
by the Sierpi\'{n}ski's characterization of 
the
rational numbers (see, for example, \cite{MR4217759}), 
 the rational Urysohn universal space is 
homeomorphic to 
the space $\qq$ of rational numbers. 
\begin{cor}\label{cor:application2}
The set of all metric $d$ such that 
$(\qq, d)$ is rational Urysohn universal space 
is meager in $\met(\qq)$. 
\end{cor}

\begin{ac}
The author would like to 
thank the referee for helpful comments 
and suggestions. 
\end{ac}


\bibliographystyle{amsplain}
\bibliography{bibtex/disco.bib}

\end{document}